\documentclass[11pt,twoside]{article} 
\usepackage{psfig,latexsym,amsfonts} 
\newtheorem{Thm}{Theorem}[section]
\newtheorem{Conj}[Thm]{Conjecture}
\newtheorem{Prop}[Thm]{Proposition}
\newtheorem{Lem}[Thm]{Lemma}
\newtheorem{Cor}[Thm]{Corollary}

\newtheorem{Def}[Thm]{Definition}
\newtheorem{Rem}[Thm]{Remark}
\newtheorem{Ex}[Thm]{Example}
\newtheorem{Exs}[Thm]{Examples}

\newcommand{\BoP}[1]{\noindent {\sc Proof#1: }}	
\newcommand{\EoP}{\hfill$\Box$\vspace{6pt}}	
\newcommand{\NoP}{\hfill$\Box$}			
\newcommand{\bb}{\mathbb} 
\newcommand{\rdrei}{{\bb R}^3}

\title{Billiard knots in a cylinder}
\author{C.\,Lamm, D.\,Obermeyer   \\   \\ 
{\footnotesize Mathematisches Institut, Beringstr.~1,} \\ 
{\footnotesize 53115 Bonn, Germany,} \\ {\footnotesize e-mail: lamm@math.uni-bonn.de}}

\begin{document}
\parindent0.5cm
\date{}

\maketitle

\vspace{1.5cm}
\begin{quotation}
\abstract{We define cylinder knots as billiard knots in a cylinder.
We present a necessary condition for cylinder knots: after 
dividing cylinder knots by possible rotational symmetries we obtain
ribbon knots. We obtain an upper bound for the number of cylinder knots
with two fixed parameters (out of three).
In addition we prove that rosette knots are cylinder knots.}
\end{quotation}

\pagestyle{myheadings}
\markboth{\sc C.\,Lamm, D.\,Obermeyer}{\sc Billiard knots in a cylinder}

\section{Introduction to billiard knots}

Billiard knots were introduced in the articles \cite{JP} 
and \cite{Lamm}. 
They are periodic billiard trajectories without self-intersections in some billiard room
in $\rdrei$. One case mentioned by V.~Jones and J.~Przytycki in \cite{JP} seems 
especially interesting to us: the case of billiard knots in a cylinder.
In this article we derive a necessary condition for cylinder knots which shows
that not all knots can be realized by them. First of all, we give a formal definition
of billiard knots.

\begin{Def} 
{\rm
Let $M$ be a 3-manifold in $\rdrei$ with piecewise smooth boundary $\partial
M$. The knot $K \subset \rdrei$ is a {\sl billiard knot in} $M$ if it is a 
polygon with

\begin{enumerate}
\item[i)]
$K \subset M$,

\item[ii)]
$\{$vertices of $K\} \subset \{p \in \partial M|
\partial M$ is smooth at $p\}$,
\item[iii)]
at the vertices $v$ we have a reflection at the tangent plane 
$T_v(\partial M)$ as known from light-rays or billiard balls.
\end{enumerate}
Also, we define a {\sl billiard link in M} as a collection of 
non-intersecting billiard knots in $M$.
}
\end{Def}

\begin{Exs}
{\rm
1.) $M={\bf D}^3$: billiard knots in a standard 3-ball. Billiard curves lie in
a plane and hence we only get trivial knots. However non-trivial links
like the Hopf links are obtainable. Every component of such a link is
trivial and the linking number of two of them is $-$1, 0 or 1. It is as yet
unknown which links exactly occur in the standard ball.

\medskip
\noindent
2.) $M={\bf I}^3$ (cube). As explained in \cite{JP} and \cite{Lamm}
billiard knots in a cube are the same as Lissajous knots defined in
\cite{Bogle}. (Recently we have found that the relationship between
Lissajous curves and billiard curves in cube was already well known
to mathematicians who study billiards, see for instance \cite{Fomenko}, p.~294,
and \cite{GZ}). The Alexander polynomial of a Lissajous knot
is a square modulo 2; from this we concluded that its Arf invariant vanishes
and that algebraic knots (e.\,g. torus knots) are not Lissajous knots \cite{Lamm}.

}
\end{Exs}

\section{Elementary properties of cylinder knots}
\begin{Def}
{\rm
A {\sl cylinder knot}  is a billiard knot in ${\bf D}^2 \times {\bf I}
=\{(z,x)\in {\bb C} \times {\bb R}:|z|\le 1,0 \le x \le 1\}$.
}
\end{Def}

We project the billiard curve on ${\bf D}^2$ and consider one edge
of the projection, see Figure \ref{diagram}.
Together with the center $M (=0 \in {\bb C})$ of ${\bf D}^2$ this edge spans a triangle; let
$\gamma$ be the angle at $M$.
Taking the edge as a part of the
billiard curve, a closed curve is achieved if and only if 
$\gamma = 2 \pi \frac{s}{n}\,(s,n \in {\bb N}$); 
here and for the rest of the article we assume gcd$(s,n)=1$ and $n \ge 2s+1$. 
The natural number $s$ and $n$ are the rotation number 
and the number of reflections of the projected curve, 
respectively.

Let $f_{1,2}:[0,1]\rightarrow {\bb C}$ parametrize the projection
proportional to arclength and let $f_{1,2}(0)=1 \in {\bb C}$ and
$f_{1,2}(\frac{1}{n})=e^{2 \pi i \frac{s}{n}} \in {\bb C}$ be the
first and second vertices. Then the vertices of the closed
billiard curve are 
$f_{1,2}(\frac{k}{n})=e^{2 \pi i \frac{s}{n} k}$, $k=0,\ldots,n-1$.

\bigskip
\centerline{\psfig{figure=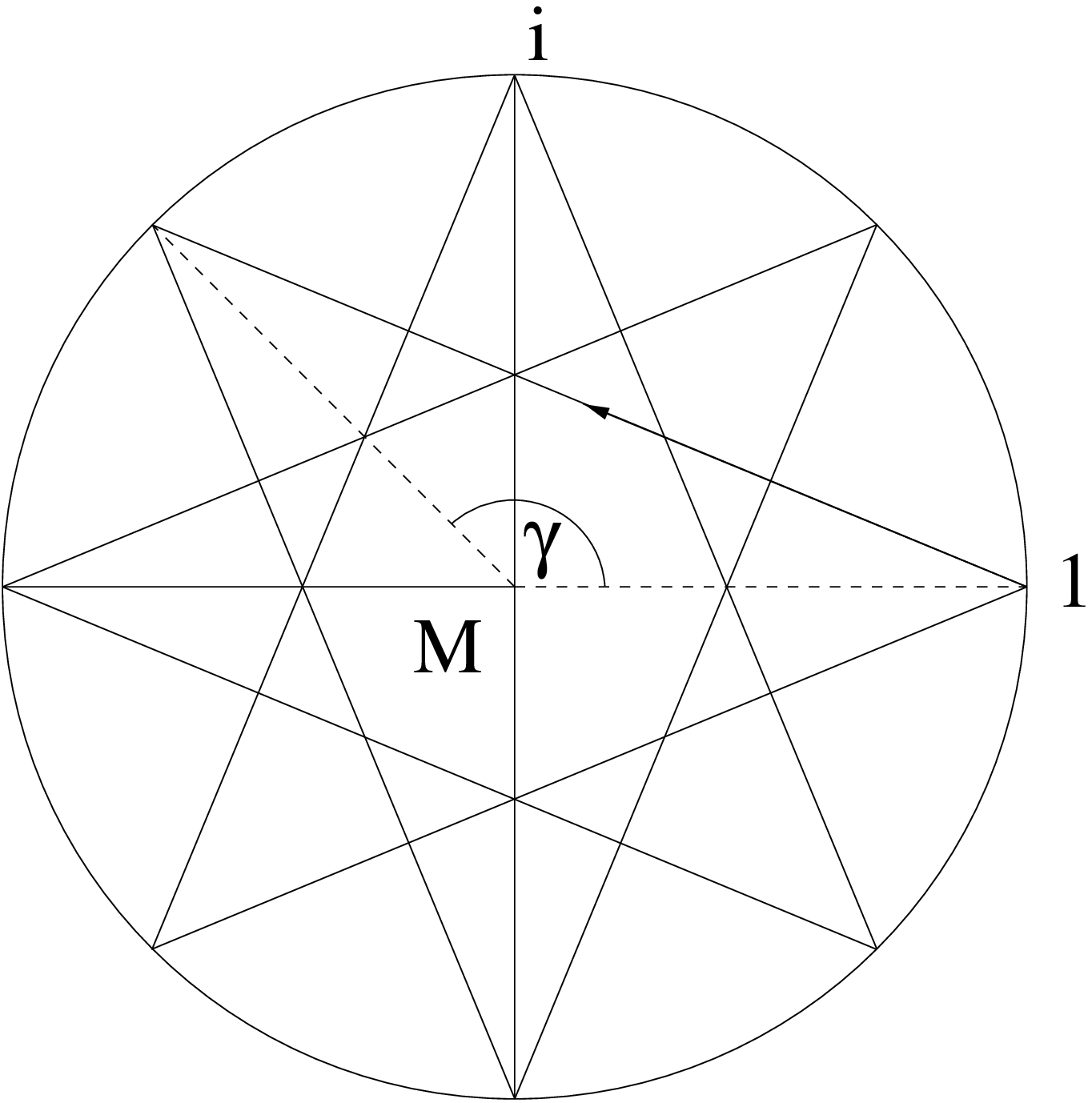,height=7cm}}
\begin{figure}[hbtp]
\caption{The projected billiard curve for $n=8$ and $s=3$.}
\label{diagram}
\end{figure}

A {\sl maximum} is a point where $K$ reflects off the ceiling.
The vertical movement is given by $m$, the number of maxima, and a 
phase $\phi$ which determines the placement of the maxima.
If $g$ is the sawtooth-function $g(t)=2|t-\lfloor t \rfloor -\frac{1}{2}|$ we
can write the height-function of the billiard curve as
$f_3 (t)=g(mt+\phi)$.
We denote the resulting billiard curve (not the knot type) by $Z(s,n,m,\phi)$.

\begin{Ex}
{\rm
If $q \ge 2p+1$ the torus knots $t(p,q)$ are cylinder knots \cite{JP}.
They are realized by $Z(p,q,q)$.
}
\end{Ex}

In the last line we neglected the phase on purpose, because 
surprisingly cylinder knots are essentially independent of it:

\begin{Lem}\label{idp_phase}
Cylinder knots are independent of the phase (up to taking mirror image).
\end{Lem}

\BoP{}
We start from a phase which gives a cylinder knot without self-intersections, 
and we show that by pushing the maxima along the knot, 
we do not change the knot up to taking mirror image. 

Let us push one maximum into the direction of the knot's orientation. The knot does not
change until, for the first time, the maximum reaches a point $A$ for which a singularity 
occurs at a crossing point $P \in D^2$. 
Let $t_1, t_2 \in [0,1]$ be the parameters of $P$; hence 
$f_{1,2}(t_1)=f_{1,2}(t_2)=P, \, f_3(t_1)=f_3(t_2)$.
By the symmetry of the function $g$ we have 
$g(mt_1+\phi) = g(mt_2+\phi) \Leftrightarrow m(t_1-t_2) \in {\bb Z}$
or $m(t_1+t_2)+2\phi \in {\bb Z}$.
If $m(t_1-t_2) \in {\bb Z}$ then the singularity cannot be removed
by changing $\phi$. If $m(t_1+t_2)+2\phi \in {\bb Z}$ then a 
sufficiently small change of $\phi$ removes the singularity.
These two cases correspond to a) and b) in Figure \ref{indep_phase}, respectively.
Now a) implies there was already a singularity at $P$ before pushing the 
maximum along the knot. Hence $P$ must look like b). We have
for all $\delta \in [0,1]$ the equation $f_3(t_1+\delta)=f_3(t_2-\delta)$. Therefore
the knot is symmetric with respect to the plane spanned by $P$ and the central
axis of ${\bf D}^2 \times {\bf I}$. By $A'$ we denote the mirror image of the
maximum at $A$. The symmetry implies that there is a maximum at $A'$ as well.
If now we continue pushing the maximum at $A$
along the knot, from $A'$ we push the other maximum into the mirror 
image of the region we came from (see Figure \ref{indep_phase}c), 
and we get the mirror image of the knot we started from. 
\EoP

\bigskip
\centerline{\psfig{figure=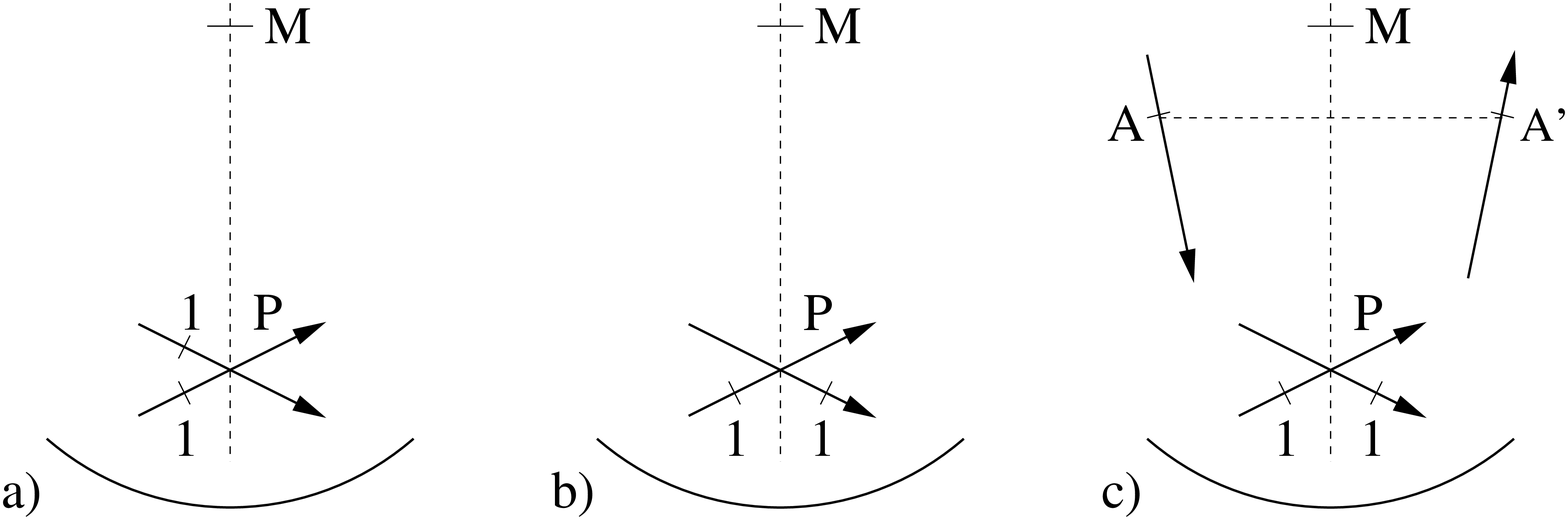,height=4cm}}
\begin{figure}[hbtp]
\caption{a) A non-removable singularity, b) a removable singularity, c) global
symmetry in case b)}
\label{indep_phase}
\end{figure}

As already mentioned the configuration of Figure \ref{indep_phase}a) 
leads to a self-intersection at $P$ for every phase. 
However we can conclude that it does not occur by using the following
lemma from number theory.

\begin{Lem}\label{independent_tan}
If $\alpha, \beta \in {\bb Q}$ with $0 < \beta < \alpha < \frac{1}{2}$ satisfy
$\tan(\pi \alpha)=\lambda \tan(\pi \beta)$ for some $\lambda \in {\bb Q}$, then
$\alpha = \frac{1}{3}, \beta=\frac{1}{6}, \lambda=3$.
\end{Lem}

\BoP{}
Let $\alpha=\frac{p_1}{q_1}$, $\beta=\frac{p_2}{q_2}$ with $\gcd(p_1,q_1)=\gcd(p_2,q_2)=1$
and $\tilde \alpha = e^{2 \pi i \alpha}$, $\tilde \beta = e^{2 \pi i \beta}$.
Then 
$
\frac{\tilde \alpha -1}{\tilde \alpha +1}=i \tan(\pi \alpha), \, 
\frac{\tilde \beta -1}{\tilde \beta +1}=i \tan(\pi \beta). 
$
Hence the cyclotomic fields ${\bb Q}[\tilde \alpha]$ and 
${\bb Q}[\tilde \beta]$ are equal.
This is possible only in the two following cases.

\noindent
First case: $q_1=q_2$. This is a contradiction to the fact that
$\cot(\pi \alpha)$ and $\cot(\pi \beta)$ are rationally independent by \cite{Wang}.

\noindent
Second case: $q_2=2 q_1$, $q_1$ odd. 
Then $-\tilde \beta$ is conjugate to $\tilde \alpha$ and 
$\frac{-\tilde \beta -1}{-\tilde \beta +1}$ is conjugate to 
$\frac{\tilde \alpha -1}{\tilde \alpha +1}$.
The degree of the extension is $\varphi(q_1)$ and we find
that the norm of $\frac{\tilde \alpha -1}{\tilde \alpha +1}$
is $\lambda^\frac{\varphi(q_1)}{2}$ 
because $\frac{-\tilde \beta -1}{-\tilde \beta +1}
\frac{\tilde \alpha -1}{\tilde \alpha +1}
=(1/\frac{\tilde \beta -1}{\tilde \beta +1})
\cdot
\frac{\tilde \alpha -1}{\tilde \alpha +1}
=\lambda$.
But we know that this norm is
either 1 or a prime number, because $\tilde \alpha - 1$ generates
a prime ideal in the ring of integers of ${\bb Q}[\tilde \alpha]$
if $q_1$ is a prime-power and it is a unit otherwise, see \cite{Wash}.
Hence $\frac{\varphi(q_1)}{2}=1$ and this implies $q_1=3$ because
3 is the only odd number $q$ with $\varphi(q)=2$. 
\EoP

\begin{Prop}\label{pathol}
In a billiard curve in the cylinder $Z(s,n,m,\phi)$ every singularity is
removable by changing $\phi$. 
\end{Prop}

\noindent
\begin{minipage}{7cm}
\BoP{} We need to show that if $t_1$ and $t_2$ are the parameters at a crossing
$P_b \,(0 < b < s)$ then there is no $m$ so that $m(t_1-t_2)\in {\bb Z}$, i.\,e. that
$t_1-t_2 \in {\bb Q}$.
We scale the projection so that its total length is 1. 
Using the triangle $MP_0P_s$ we get $|MP_0|=1/(2n\tan(\frac{\pi s}{n}))$ because
$|P_0P_s|=\frac{1}{2n}$, and from the triangle $MP_0P_b$ we compute 
$|P_0P_b|=\frac{1}{2n}\frac{\tan(\frac{\pi b}{n})}{\tan(\frac{\pi s}{n})}$. 
\end{minipage}
\hfill
\begin{minipage}{7cm} 
\centerline{\psfig{figure=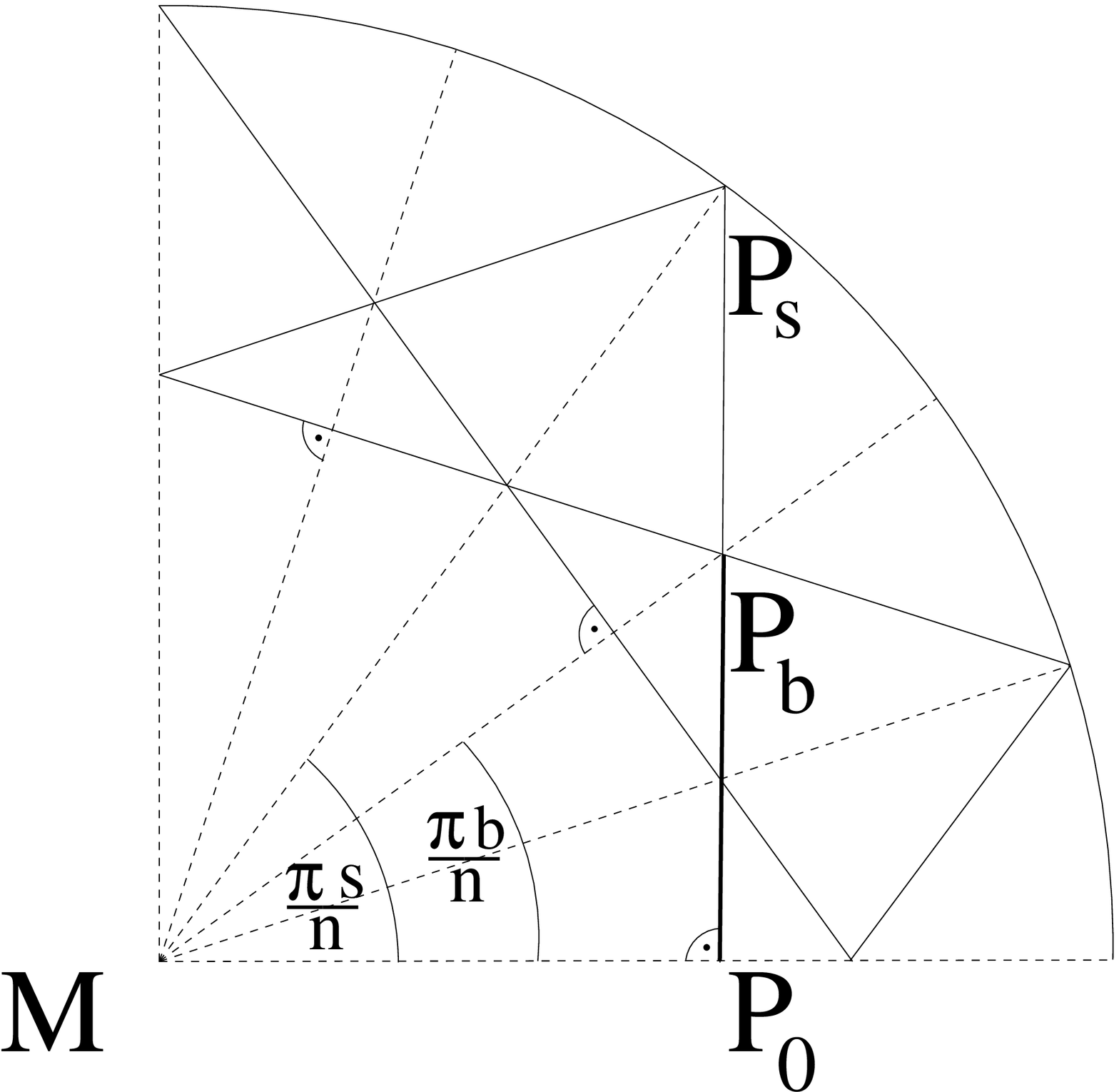,height=5cm}}
\end{minipage} 

\medskip
So there are $k_1,k_2\in {\bb Z}$ with
$$
t_1=\frac{k_1}{2n}+\frac{1}{2n}\frac{\tan(\frac{\pi b}{n})}{\tan(\frac{\pi s}{n})},\,\,
t_2=\frac{k_2}{2n}-\frac{1}{2n}\frac{\tan(\frac{\pi b}{n})}{\tan(\frac{\pi s}{n})}.
$$
From this it is immediate that if $m(t_1-t_2)\in {\bb Z}$ then 
$$
\tan (\frac{\pi b}{n})/\tan (\frac{\pi s}{n}) \in {\bb Q}.
$$
Lemma \ref{independent_tan} shows that the only
solution is $n=6b, s=2b$ which contradicts gcd($s,n) = 1$.
\EoP

\medskip
From now on we will denote a cylinder knot by $Z(s,n,m)$ but keep in
mind that the knot-type is determined by this notation only up to mirror image.
However, the billiard {\sl curve} is denoted by $Z(s,n,m,\phi)$.

We take a look at the symmetries of cylinder knots
(compare \cite{JP}, Theorem 3.11).
A knot $K$ has {\sl cyclic period n} if $K$ is fixed under a rotation
of $\frac{2\pi}{n}$ around an axis $h$ with
$h\cap K = \emptyset$. The {\sl linking number} associated to this
rotation is $lk(h,K)$. A knot $K$ is {\sl strongly positive amphicheiral}
if there is an involution $i: \rdrei \rightarrow \rdrei$ with
$i(K) = K$, preserving the orientation of the knot and reversing
the orientation of $\rdrei$.

\begin{Lem}\label{strong}
a) If $gcd(n,m) = d > 1$, then $Z(s,n,m)$ has cyclic period $d$ with linking
number $s$. Conversely, if $Z(s,n,m, \phi)$ has rotational symmetry $d'$ then
$d' | d$.

\noindent
b) If $n$ is even and $m$ is odd, then $Z(s,n,m)$ is strongly positive
amphicheiral.
\end{Lem}

\BoP{}
a) A rotation of $\frac{2 \pi}{n}$ about the central axis corresponds to
$t \mapsto t+\frac{k}{n}$ in the parametrization, where $ks \equiv 1$ (mod $n$).
Hence for $d|n$ the shift $t \mapsto t+\frac{k}{d}$ gives a rotation of $\frac{2 \pi}{d}$,
which maps the projection to itself. But because $d|m$ we have $f_3(t)=f_3(t+\frac{k}{d})$;
this shows that the curve $Z(s,n,m,\phi)$ is rotationally symmetric with period $d$.
For the converse let $d' \in {\bb N}$ be a rotational symmetry of the
diagram. Then for the projection we need $d'|n$ and for the height function
$\frac{1}{d'} \in \frac{1}{m}{\bb Z}$. Then $d'|m$ and hence $d'|d$.

In the case of b) a rotation of $\pi$ preserves the diagram but exchanges maxima
and minima: the curve is symmetric with respect to $(0,\frac{1}{2})$, the center of
the cylinder.
\EoP

The number $d$ is the maximal period of the diagram. 
We call it the {\sl maximal billiard period}.

\begin{Rem}\label{factor_sym}
{\rm
In the periodic case the {\sl factor knot} $K^{(d)}$ is obtained by identifying the
radial faces of a $\frac{2\pi}{d}$-slice of the cylinder. We do this
by a process of enlarging the slice to the whole cylinder, that is (in
cylinder coordinates) mapping the angle $\alpha$ to $\alpha \cdot d$.
The result is no longer a billiard knot since its lines are not straight
but it still has those symmetries of cylinder knots which are derived
from the dihedral symmetries of their projections.
}
\end{Rem}

\section{A necessary condition for cylinder knots}
In this section we prove that not every knot can be a cylinder 
knot. In order to show this, we will prove that after factoring by
the maximal billiard period, the resulting factor knot is a ribbon
knot.

We want to allow maxima and minima on the boundary of ${\bf D}^2$
as well, although our definition of billiard knots forbids this because
the boundary of ${\bf D}^2\times {\bf I}$ is not smooth there. Instead we
start from the projected curve and distribute the maxima equidistantly.

\begin{Lem}\label{singsym}
a) A maximum or minimum at a vertex $P$ on the boundary of ${\bf D}^2$ yields
singularities at every crossing point on the plane spanned by $P$ and the
central axis of ${\bf D}^2 \times{\bf I}$.

\noindent
b) If $\gcd(n,m)=1$ then there are no singular crossings other than the
singularities mentioned in a).
\end{Lem}

\BoP{}
a) The knot projection on ${\bf D}^2$ is symmetric with respect to the
plane spanned by $P$ and the central axis $0 \times{\bf I}$. Now if we have a maximum or
minimum at $P$ then the height function is symmetric as well, so the whole
billiard curve is symmetric with respect to the plane. This can only happen
if the crossing points on the plane are self-intersections.

\noindent
b) Similarly to a), a singularity implies (Figure \ref{indep_phase}, Proposition
\ref{pathol}) that the billiard curve is symmetric
with respect to the plane spanned by the crossing point and the central
axis. Hence a singular crossing outside the symmetry plane of a) yields
a second symmetry plane of the billiard curve. The product of these two
reflections is a rotation around the central axis fixing the curve. But
since $\gcd (n,m)=1$ this is impossible (Lemma \ref{strong}).
\EoP{}

A \emph{ribbon knot} is a knot which is the boundary of a singular disk with
only ribbon singularities. 

\bigskip
\centerline{\psfig{figure=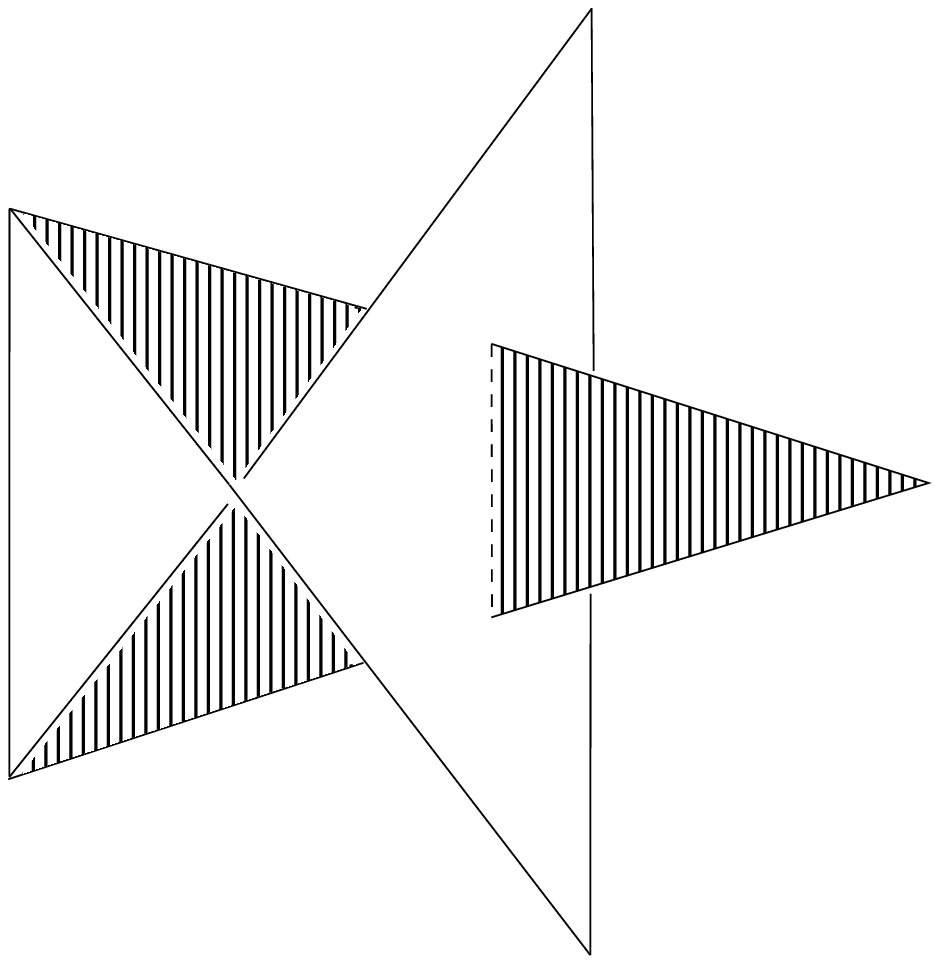,height=7cm}}
\begin{figure}[hbtp]
\caption{A ribbon for a cylinder knot diagram}
\label{example}
\end{figure}


\begin{Prop}\label{ribbon}
Let $K=Z(s,n,m)$ be a cylinder knot and $d=\gcd(n,m)$ its maximal billiard period.
Then the factor knot $K^{(d)}$ is a ribbon knot.
\end{Prop}

\BoP{}
Remember that $K^{(d)}$ still admits the symmetries of cylinder knots
(Remark \ref{factor_sym}). The idea to prove this proposition
is to connect points with their mirror
image as $A$, $A'$ in Figure \ref{indep_phase}c). This is similar to the well-known
technique for $K \sharp (-K^*)$.
 
Let $f_{1,2}: [0,1] \rightarrow {\bf D}^2 \subset {\bb C}$ be the mapping describing the 
projection of $K^{(d)}$ with $f_{1,2}(0)=f_{1,2}(1)=1$, 
and $f_3: [0,1] \rightarrow {\bf I} $ the
height function starting with a maximum at 1: $f_3(1)=1$. 
By the preceding lemma
we get singularities on ${\bb R} \times {\bf I}$ and there are no further
singularities. 

Now by the following steps we construct a ribbon for $K^{(d)}$:
In step 4) the singularities on ${\bb R} \times {\bf I}$ are 
removed so that they unfold to half-twists in the ribbon.

\bigskip
1) We start with a mapping
\begin{eqnarray*} 
B : [0,1]^2 &\rightarrow& {\bf D}^2 \times {\bf I} \\
(x,y) &\mapsto & ( (1-y)f_{1,2}(\frac{x}{2}) + 
y\overline{{f}_{1,2}(\frac{x}{2})},
f_3(\frac{x}{2}) ).
\end{eqnarray*}

Let $S=$\,im$(B)$ be the singular disk described by $B$.
$B(x,0)$ parametrizes the first half of $K^{(d)}$; the second half is 
parametrized by $B(x,1)$ since $K^{(d)}$ is symmetric with respect to 
${\bb R} \times {\bf I}$. 
$B(0,y)$ and $B(1,y)$ are constant. Hence $K^{(d)}$ is the boundary of $S$.

2) Because of the symmetry of $B$, all the self-intersections of $S$ 
are ribbon singularities.

3) If three or more parts of the ribbon meet in one arc, we may deform $S$
in a sufficiently small neighbourhood to push away
one of the parts without destroying the symmetry.

4) Now let us finally push away the maximum from 1 to remove the singularities
on ${\bb R} \times {\bf I}$. One of the two sides
will yield the desired knot (the other its mirror image, see the proof of
Lemma \ref{idp_phase}). As long as the
maximum is sufficiently close to 1, we do not change
the knot outside ${\bb R} \times {\bf I}$ (by Lemma \ref{singsym}, b)), 
and we move $S$
together with the knot. The singularities on ${\bb R} \times 
{\bf I}$ open to form half-twists in $S$, and the result is
the desired ribbon for our factor knot (Figure \ref{example}).
\EoP{}

\bigskip

\centerline{\psfig{figure=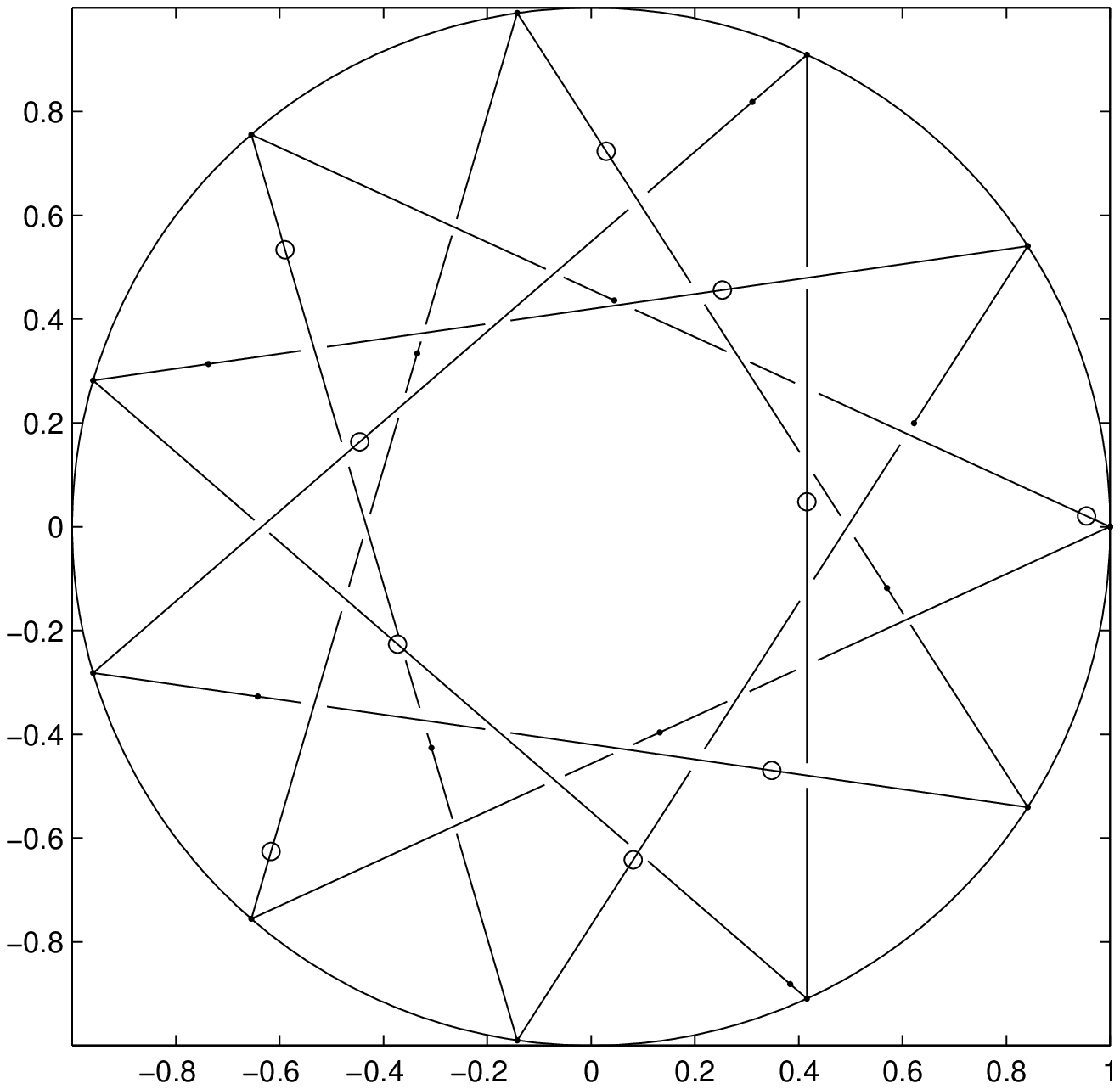,height=10cm}}
\begin{figure}[hbtp]
\caption{An example of a ribbon cylinder knot: $Z(4,11,10)$. The circles and dots
are the positions of the maxima and minima.}
\label{exeins}
\end{figure}
 
\noindent
{\bf Remark:} We observe that $K^{(d)}$ is a generalized union of
a knot diagram and its mirror image as defined in \cite{KT} and \cite{L2}.

\bigskip
We sum up what we know about cylinder knots:

\begin{Thm}\label{condition}
Let $K=Z(s,n,m)$ be a cylinder knot and $d=$ gcd$(n,m)$ ($d=1$ may happen).
Then $K$ has cyclic period $d$ with linking number $s$ and the factor
knot $K^{(d)}$ is ribbon.  
In particular, $K^{(d)}$ satisfies:

a) The determinant is a square,

b) the Arf invariant vanishes,

c) the signature vanishes.

\noindent
If $\frac{n}{d}$ is even then $K^{(d)}$ is also strongly positive amphicheiral. 
\end{Thm}

\BoP{} 
By Proposition \ref{ribbon} $K^{(d)}$ is a ribbon knot.
We know from \cite{Kaw} that its
signature vanishes and the Alexander polynomial has the form \[
\Delta_{K^{(d)}}(t) = F(t)F(t^{-1}) \] for some $F \in \,{\bb Z}[t]$.
Since det$(K^{(d)})=\Delta_{K^{(d)}}(-1)$ we have a)
and c). Part b) is a corollary from a) because 1 is the only odd square 
modulo 8, hence det$(K^{(d)})\equiv \pm1$ (mod 8) and the Arf invariant vanishes.

The diagram of $K^{(d)}$ has ${\frac{n}{d}}$ vertices and ${\frac{m}{d}}$
maxima. The numbers ${\frac{n}{d}}$ and ${\frac{m}{d}}$ are coprime;
in particular, if ${\frac{n}{d}}$ is even then ${\frac{m}{d}}$ has to be odd,
hence by Lemma \ref{strong} $K^{(d)}$ is strongly positive amphicheiral.
\EoP{}

\begin{Cor}
Cylinder knots are either periodic or ribbon. Hence not all knots 
are cylinder knots.
\end{Cor}

For instance, the knot $8_{10}$ (notation as in \cite{Kaw}) has no 
period and its determinant is 27, hence it does not match our necessary
condition. 
But not even all periodic knots are cylinder knots. The knot $5_2$ has
2 as its only period, but the linking number with its axis of rotation
is $\pm$1; furthermore, det($5_2$)=7, so $5_2$ cannot be a cylinder knot.
Note that $5_2$ is the first nontrivial Lissajous knot. In fact, we know
only one example of a Lissajous knot which is at the same time a cylinder 
knot: this is $3_1\sharp3_1^* = Z(3,11,4)$.

The argument of the last paragraph shows:

\begin{Cor}
Let $K$ be a knot which is not ribbon.
If the complete list of its cyclic periods together with the
linking numbers with the respective axis of rotation is
$((q_1,\lambda_1),\ldots,(q_k,\lambda_k))$ and if for all $i$  we have
$|\lambda_i| < br(K)$ then $K$ is not a cylinder knot.\NoP
\end{Cor}

\section{An upper bound for the number of cylinder knots with given
$s$ and $n$}
If $Z(s,n,m)$ is a cylinder knot we write 
$$
Z(s,n,m)=\prod_{i=1}^n \Bigl[\prod_{j \,{\rm odd} \atop j < s} 
\sigma_j^{\varepsilon_{i,j}}
\prod_{j \,{\rm even} \atop j < s} 
\sigma_j^{\varepsilon_{i,j}}\Bigr]~^{\bf \widehat{}}\,,\,\,\varepsilon_{i,j} \in \{\pm 1\},
$$
reading the braid in counterclockwise direction. Of course
there is a choice of rotating and mirroring the braid, but for
Theorem \ref{cycl} we fix one of the braids.
We consider the first subscript $i$ of $\varepsilon_{i,j}$ (and $x_{i,j}$,
$y_{i,j}$ in the proof of Theorem \ref{cycl}) to be in ${\bb Z}/n{\bb Z}$.

If $n$ and $m$ are coprime, Theorem \ref{cycl} says that for 
each $j$ we have a series of positive and negative crossings.
These series have the same length (for even $n$), or their lengths differ only by one.

\begin{Thm}\label{cycl}
Let $\gcd(n,m)=1$ and $a\in\{1,\ldots,n-1\}$
with $am \equiv s \pmod{n}$. Then for each $j\in\{1,\ldots,s-1\}$ there
is exactly one $b_j\in\{1,\ldots,n\}$ so that for $n$ even
$$\varepsilon_{b_j+i\cdot a,j} = \Bigl\{
\begin{tabular}{rl}
$+1$ & for $i=0,\ldots,\frac{n}{2}-1$, \\
$-1$ & for $i=\frac{n}{2},\ldots,n-1$,
\end{tabular}
\Bigr.
$$
\mbox{ and for n odd}

\vspace{-0.6cm}
$$
\varepsilon_{b_j+i\cdot a,j} = \Bigl\{
\begin{tabular}{rl}
$+1$ & for $i=0,\ldots,\frac{n-1}{2}$, \\
$-1$ & for $i=\frac{n+1}{2},\ldots,n-1$,
\end{tabular}
\Bigr.
\,\,{or }\,\,
\Bigl\{
\begin{tabular}{rl}
$-1$ & for $i=0,\ldots,\frac{n-1}{2}$, \\
$+1$ & for $i=\frac{n+1}{2},\ldots,n-1$.
\end{tabular}
\Bigr.
$$
\end{Thm}

\BoP{}
We denote the crossing corresponding to $\sigma_j^{\varepsilon_{i,j}}$ 
by $(i,j)$ and the parameters of the crossing $(i,j)$ by
$x_{i,j}$ and $y_{i,j}$. We assume that the $x_{i,j}$ belong to the left string
at a crossing and the $y_{i,j}$ to the right string, respectively. 
For the proof of the theorem we fix one $j$ and omit it from the subscripts for simplification.
Let $k \in \{1,\ldots,n-1\}$ as in the proof of Lemma \ref{strong} satisfy
$ks \equiv 1$ (mod $n$).
Because of $am \equiv s$ (mod $n$) we get the following congruences.
\begin{eqnarray*}
x_{i+1} &\equiv& x_i + \frac{k}{n} \pmod{1},\\
x_{i+a} &\equiv& x_i + \frac{ak}{n} \pmod{1},\\
mx_{i+a} &\equiv& mx_i + \frac{amk}{n}\equiv mx_i + \frac{1}{n} \pmod{1}.
\end{eqnarray*}
In the same way we get $my_{i+a} \equiv my_i + \frac{1}{n} \,(\mathrm{mod}\,\,1)$. The height
function is $f_3(t)= g(mt+\phi)$. Since the multiplication by $m$ is already done in
the computation above, the heights of the $mx_{1+l\cdot a}$ and of the $my_{1+l\cdot a}$ for
$l=0,\ldots,n-1$ differ only by a constant, see the figure. 

\bigskip

\bigskip

\centerline{\psfig{figure=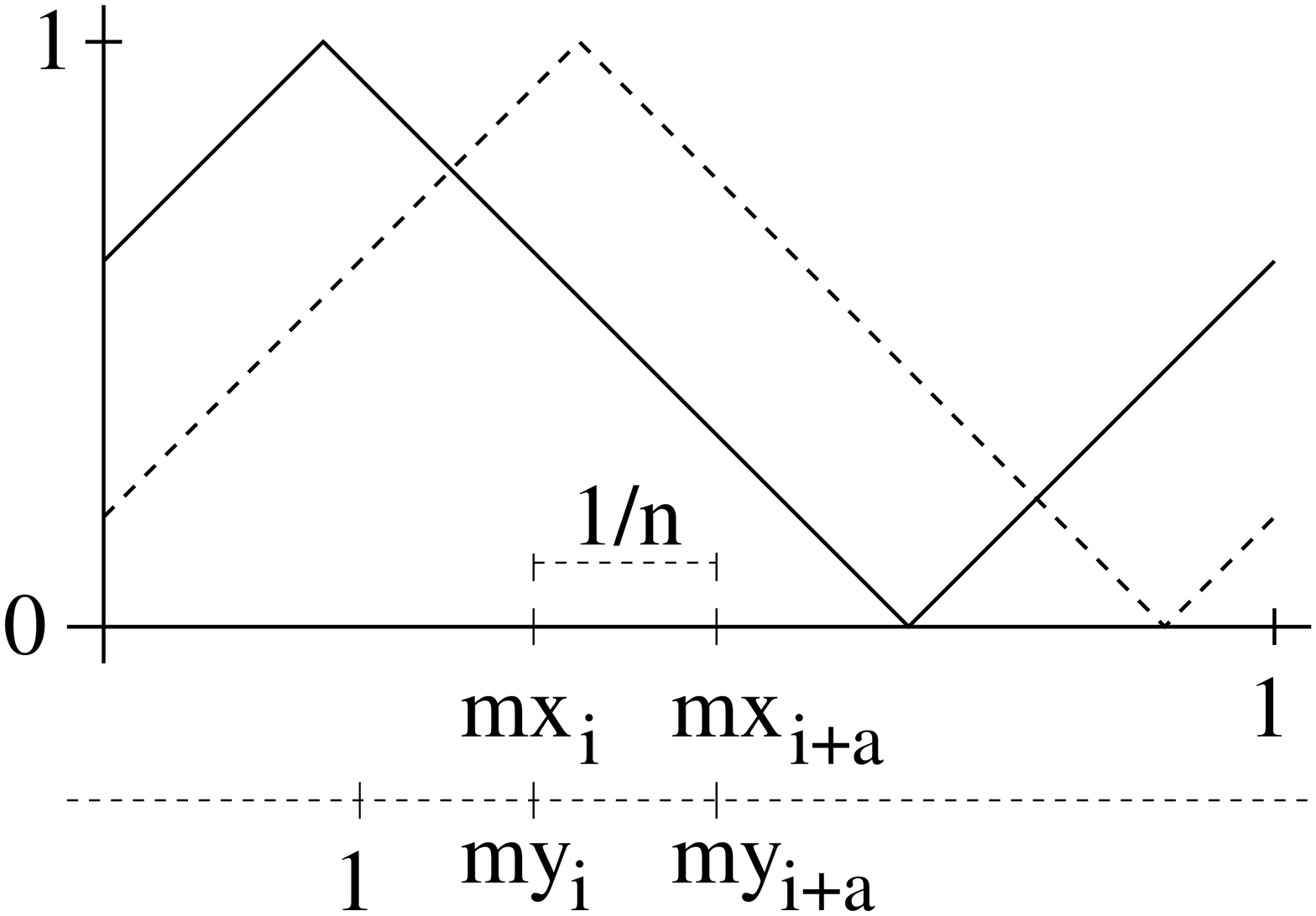,height=4cm}} 

\bigskip

For even $n$ we choose for $b_j$ the uniquely determined index for which
$g(mx_{b_j}+\phi) > g(my_{b_j}+\phi)$ and $g(mx_{b_j-a}+\phi) < g(my_{b_j-a}+\phi)$.
For odd $n$ we choose the $b_j$ which is at the beginning of the longer series of
identical signs (these can be positive or negative, as stated in the theorem).
\EoP

\noindent
Notation: We define $v_j:=(\varepsilon_{1,j},\ldots,\varepsilon_{n,j}) \in \{+1,-1\}^n$
for all $j \in \{1,\ldots,s-1\}$.

\noindent
In the following two theorems we identify rotated and mirrored braids.

\begin{Thm}\label{poss}
Let $\gcd(n,m)=1$. Then, given $v_1 \in \{+1,-1\}^n$, for each $j\ge 2$ there are at most two possibilities
for $v_j$, and they differ only by sign.
\end{Thm}

\BoP{}
Assume a phase $\phi$ gives the crossing-sign vector $v_j$
and the height diagram looks like the following (different cases
for even and odd $n$).

\bigskip

\bigskip

\centerline{\psfig{figure=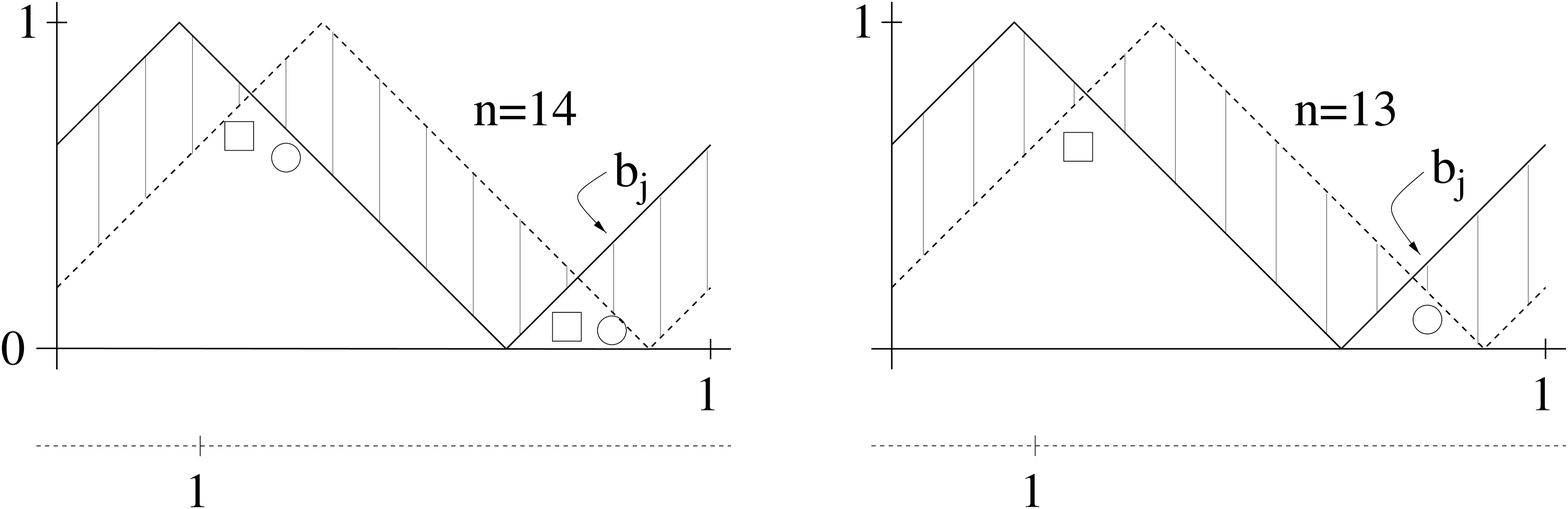,height=4cm}}

\bigskip

When we change the phase $\phi$ there are some crossings in danger
to become singular and to be changed.
\begin{enumerate}
\item
$n$ is even: Then, decreasing $\phi$ changes the pair of crossings
$\{(b_j,j),(b_j+\frac{n}{2}a,j)\}$ and increasing $\phi$ changes
the pair $\{(b_j-a,j),(b_j-a+\frac{n}{2}a,j)\}$. These two pairs are
marked by circles and squares in the above figure. We call them the
{\sl critical pairs} (of crossings) of $v_j$.
\item
$n$ is odd: In this case there is only the crossing $(b_j,j)$
which is changed for decreasing $\phi$, and for increasing 
$\phi$ the crossing $(b_j+\frac{n-1}{2}a,j)$. We call them the
{\sl critical crossings} of $v_j$.
\end{enumerate}

\bigskip
{\bf Claim:} The critical pairs (or critical crossings) of a crossing-sign
vector determine it up to sign.

\BoP{ of claim} The critical pairs or critical crossings
determine the beginning and end of the two
series of identical crossing-signs because $n \ge 5$ if $s \ge 2$. 
If $n$ is even, the series have equal
length $\frac{n}{2}$. If $n$ is odd the critical crossings belong to
the longer series (of length $\frac{n+1}{2}$). In both cases there
is only the choice of a global sign.\,$\Box$

This fact yields the proposition for odd $j$ because the singular phases for 
$v_1$ coincide with those for $v_j$. 

For odd $n$ and even $j$ a singularity at $(i,1)$ implies a singularity
at $(i+\frac{n-1}{2},j)$, hence the critical crossings determine
each other.

For even $n$ and even $j$ we use a different method: we
place a maximum in the middle of a chord at $P$. Then exactly two crossings
in $v_1$ are singular and the singular knot has a symmetry plane
spanned by $P$ and the central axis. 
In $v_j$ there are exactly two indices $i$ with $\varepsilon_{i,j}=-\varepsilon_{i+a,j}$.
Hence for $v_j$ there are only two possibilities $v_j'$ and $v_j''=-v_j'$. Resolving
the singularities by slightly moving the maximum away from
$P$ gives two possibilities in each case but they are mirror images of each
other.
\EoP

\medskip
\centerline{\hspace{1cm}\psfig{figure=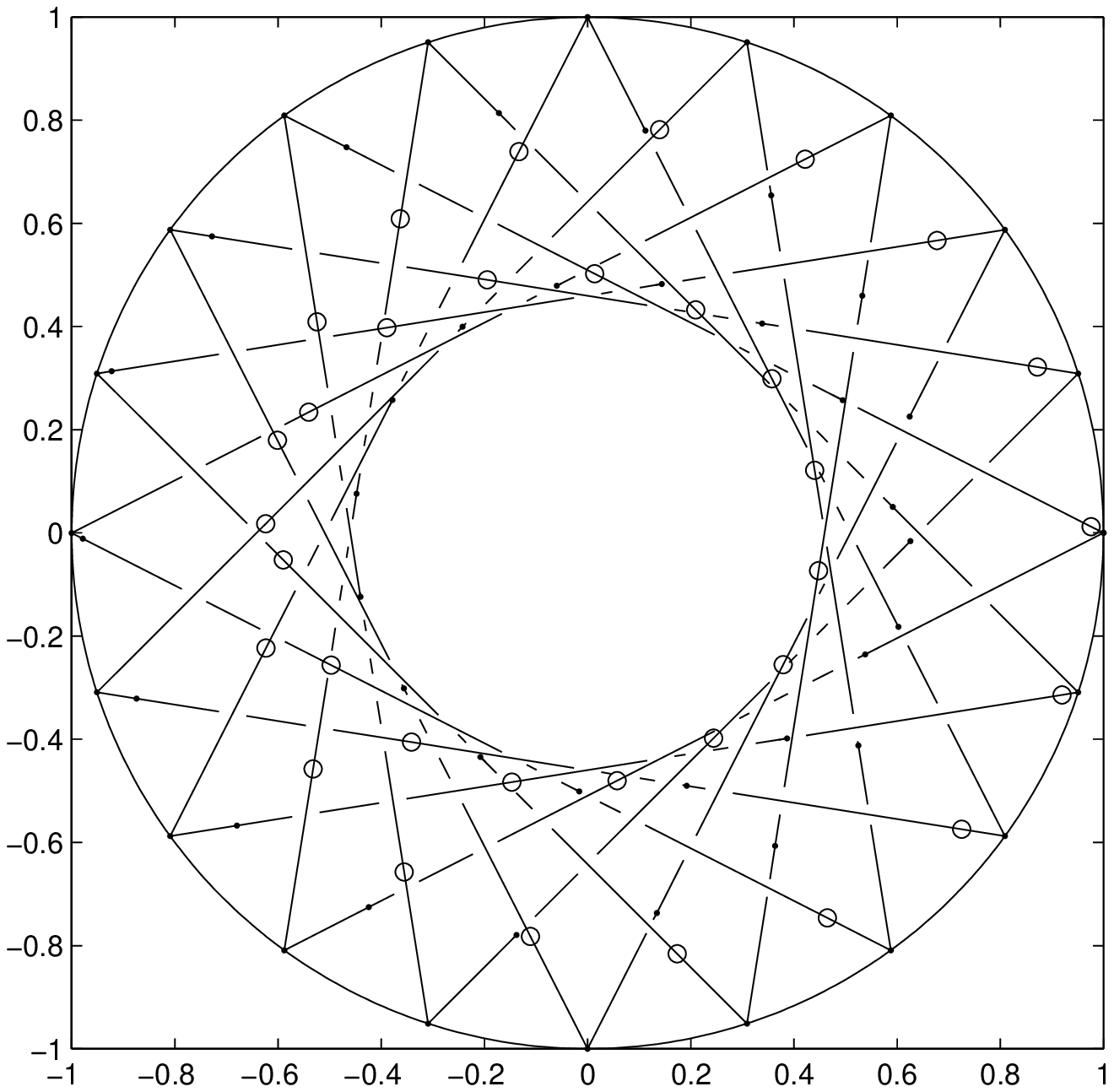,height=9cm}}
\begin{figure}[ht]
\caption{The cylinder knot $Z(7,20,33)$. It is ribbon and strongly positive 
amphicheiral.}
\label{zwanz}
\end{figure}

\begin{Thm}
For given $s, n$ there are at most $(n+1)2^{s-3}$ cylinder braids.
\end{Thm}

\BoP{}
We count all possibilities for the $v_j$ for the different 
periods of the diagram. 
Let $d=\gcd(n,m)$ be the period of the diagram and $c=\frac{n}{d}$. Then
for $c=1$ there are $2^{s-2}$ possibilities because we can choose
$\varepsilon_{1,1}=1$ and there are two possibilities for every other
$j$.
(Following Remark \ref{factor_sym}, we can apply the Theorems \ref{cycl} and
\ref{poss} to the factor knots. The only exceptional case is $n/d=4$, 
because then every crossing belongs to a critical pair. 
We check this case directly by the
method of placing a maximum in the middle of a chord or on a vertex.)
For $c\ge 2$ there are $\frac{\varphi(c)}{2}$ possibilities to choose
$a$, because $a$ and $-a$ give the same $v_1$. So there are at most
$\frac{\varphi(c)}{2}\cdot 2^{s-2}$ possibilities for all $j$. Some
care is needed in the case $c=2$. Here a braid with $v_1=v_2=(+1,-1)$ 
corresponds by a rotation around an axis in the diagram plane 
to a braid with $v_1=(+1,-1)$ and $v_2=(-1,+1)$. Hence there are 
$2^{s-3}=\frac{\varphi(2)}{2}\cdot 2^{s-2}$ possibilities.
The sum over all $c$ is
$$
2^{s-2}+\sum_{c|n \atop {c\ge 2}} \varphi(c) 2^{s-3}=2^{s-2}+(n-1)2^{s-3}=(n+1)2^{s-3}.
$$
\EoP

\section{Rosette knots}
Rosette knots have been defined by Kr\"otenheerdt \cite{Kr} as a
special class of alternating knots. These knots and the corresponding
links can all be shown to be cylinder knots (or links), so we have
quite a general class of examples.

\begin{Def}
{\rm
Let $k$ and $s$ be natural numbers. The closure of the $s$-string braid \[
(\sigma_1\sigma_2^{-1}\sigma_3\sigma_4^{-1}\ldots\sigma_{s-1}^{(-1)^{s}})^k \]
is a link with gcd($s,k$) components; we call it the {\sl rosette link} 
$R^k_s$. In particular, if $k$ and $s$ are coprime, we call $R^k_s$ a
{\sl rosette knot}.
}
\end{Def}

\begin{Rem}
{\rm
The rosette knots in the knot tables up to 10 crossings 
are the torus knots $t(2,2n+1)=R^{2n+1}_2$ for $n=1,2,3,4$, the 
figure-eight knot $4_1=R^2_3$,
$8_{18}=R^4_3$, $10_{123}=R^5_3$ and $9_{40}=R^3_4$.
For rosette links we have the well-known examples $R^2_2$ (Hopf link) 
and $R^3_3$ (Borromean rings).
}
\end{Rem}

\begin{Thm}
For $k\ge2$ and $\gcd(k,s)=1$, the rosette knot $R^k_s$ is the 
cylinder knot $Z(s,k(s+1),k)$.
\end{Thm}

\BoP{}
We start with the singular cylinder knot $K=Z(s,k(s+1),k,\phi)$ by placing a maximum
on a vertex $P_1$. Because the billiard curve has period $k$ it is
enough to look at the factor knot; we consider the slice with angle $\frac{2 \pi}{k}$ and
deform the braid to a grid pattern for better visualization. 
We divide the height interval $[0,1]$ into layers of height $\frac{2}{s+1}$
(and for even $s$ one bottom layer of height $\frac{1}{s+1}$). Then on the
symmetry planes there are singular crossings and the layering shows that
the crossings in the first half of the slice are negative and in the second
half they are positive (the convention is, that positive crossings correspond
to the $\sigma_i^{-1}$).

See Figure \ref{rosette2} for $s=9$ as an example. Here we write $A$, $B$,$\ldots$,
$E$ for the layers $[\frac{8}{10},1]$, $[\frac{6}{10},\frac{8}{10}]$,\ldots,$[0,\frac{2}{10}]$.
For instance the chord $P_1P_{10}$ belongs to the layer $A$, $O_{10}P_9$ to $B$, hence
$P_1P_{10}$ is above $O_{10}P_9$.

\medskip
\centerline{\psfig{figure=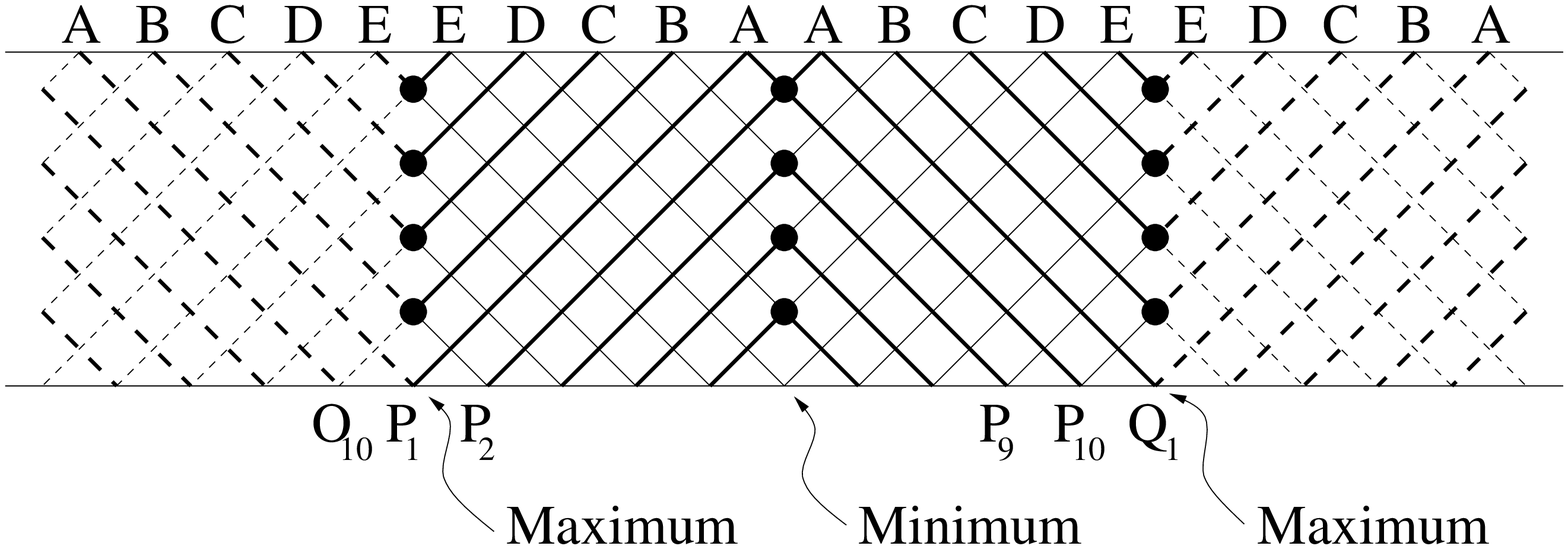,height=4.5cm}}
\begin{figure}[ht]
\caption{The braid (\ref{braid1}) for $s=9$ with singularities on the symmetry planes}
\label{rosette2}
\end{figure}

The next step is to push the maxima slightly to the left. This resolves the singularities
to positive crossings at the planes at $P_1$ and $Q_1$ and to negative crossings in the
central plane of the slice. Hence the factor knot has as a braid

\begin{equation}\label{braid1}
(\sigma_1\sigma_3\sigma_5\ldots\sigma_{s-2}\sigma_2\sigma_4\ldots
\sigma_{s-1})^{\frac{s+1}{2}}(\sigma_1^{-1}\sigma_3^{-1}\ldots\sigma_{s-2}^{-1}
\sigma_2^{-1}\sigma_4^{-1}\ldots\sigma_{s-1}^{-1})^{\frac{s+1}{2}} 
\end{equation}
if $s$ is odd, and
\begin{equation}\label{braid2}
(\sigma_1\sigma_3\sigma_5\ldots\sigma_{s-1}\sigma_2\sigma_4\ldots
\sigma_{s-2})^{\frac{s}{2}}(\sigma_1\sigma_3\ldots\sigma_{s-1}\sigma_2^{-1}
\sigma_4^{-1}\ldots\sigma_{s-2}^{-1})
(\sigma_1^{-1}\sigma_3^{-1}\ldots\sigma_{s-1}^{-1}
\sigma_2^{-1}\sigma_4^{-1}\ldots\sigma_{s-2}^{-1})^{\frac{s}{2}} 
\end{equation}
if $s$ is even.

\centerline{\hspace{1cm}\psfig{figure=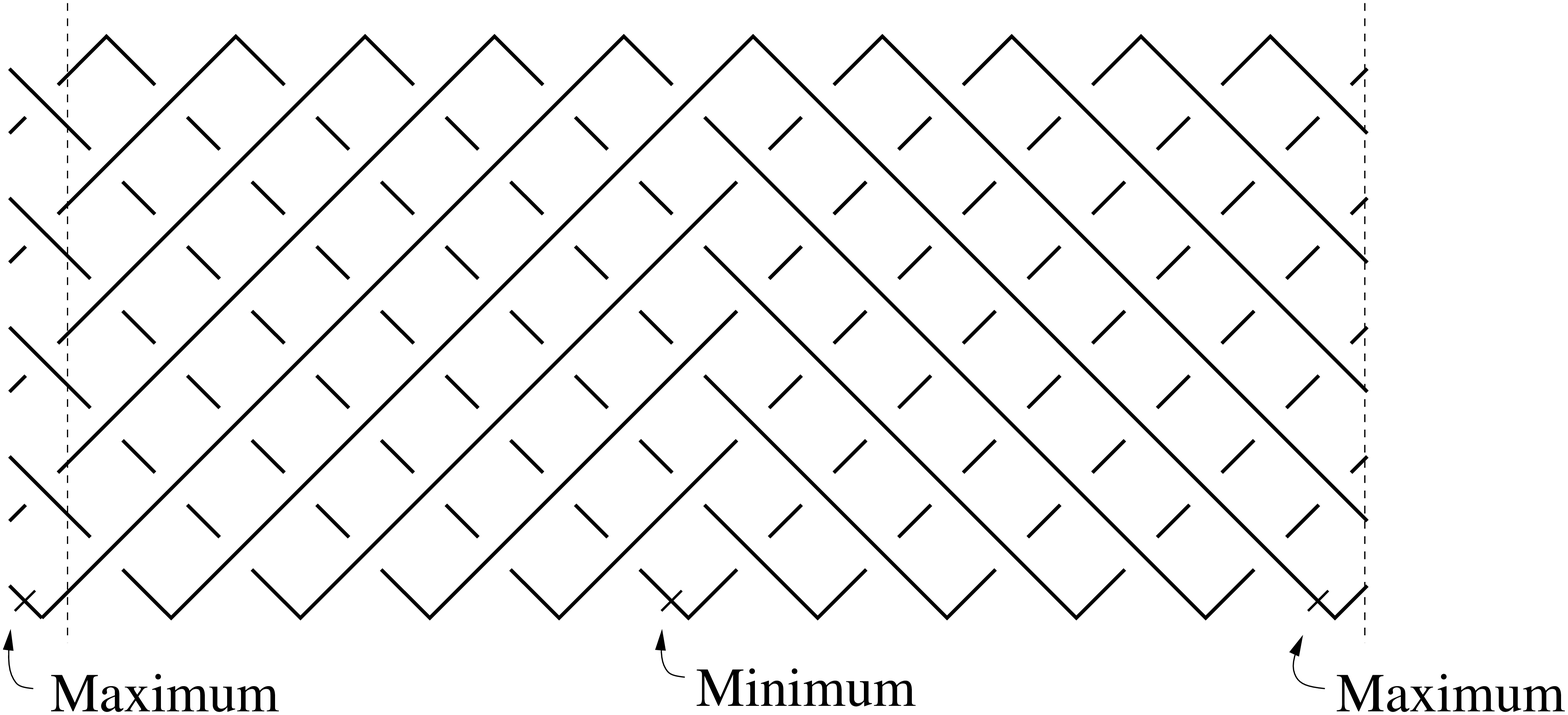,height=4.5cm}}
\begin{figure}[ht] 
\caption{The braid (\ref{braid1}) for $s=9$} 
\label{rosette3} 
\end{figure} 

We show by induction on $s$ that (\ref{braid1}) is exactly the braid 
$$
\sigma_1\sigma_3\sigma_5\ldots\sigma_{s-2}
\sigma_2^{-1}\sigma_4^{-1}\ldots\sigma_{s-1}^{-1}. 
$$

The induction step works as follows: In the braid (\ref{braid1}), 
we observe that there is one string running above all the others and one
running beneath, see Figure \ref{rosette3}.
Hence we may pull out these two which leave only one negative
crossing at the beginning and one positive at the end; the part in between then
matches our induction hypothesis if we rotate it around the horizontal axis
$P_1Q_1$ with an angle of $\pi$, and the proposition follows.
By a similar argument we see that (\ref{braid2}) is 
$$
\sigma_1\sigma_3\sigma_5\ldots\sigma_{s-1}
\sigma_2^{-1}\sigma_4^{-1}\ldots\sigma_{s-2}^{-1}\,.
$$
The reader is invited to prove that both braids are conjugate to
\[
\sigma_1\sigma_2^{-1}\sigma_3\sigma_4^{-1}\ldots\sigma_{s-1}^{(-1)^{s}}\,. \]

In the whole cylinder, we have $k$ times this braid; hence the proof is
complete.
\EoP{}

\begin{Rem}
{\rm
If gcd$(k,s)>1$, then the construction yields the desired billiard link
in the cylinder, but we have to assume that the vertices in the projection
are equidistant, and for every component the maxima are at the right 
positions; in this case we cannot apply Lemma \ref{idp_phase}.
}
\end{Rem}

From \cite{Kr} and \cite{Mura} we already know that $R^k_s$ is amphicheiral
if and only if $s$ is odd. 
Together with Lemma \ref{strong} we sharpen this condition:

\begin{Cor}
If both k and s are odd then $R^k_s$ is strongly positive amphicheiral.
\NoP
\end{Cor}
\medskip

\begin{Conj}
Any generalized rosette knot, i.\,e. the closure of \[
(\sigma_1^{\varepsilon_1}\sigma_2^{\varepsilon_2}\ldots
\sigma_{s-1}^{\varepsilon_{s-1}})^k\,, \]
with $\varepsilon_i \in \{\pm1\}$, is a cylinder knot.
\end{Conj}

These knots are periodic with period $k$, linking number $\pm s$ and 
trivial factor knots. Hence they match the conditions of Theorem
\ref{condition}. For example all torus knots are generalized rosette knots.

We mention that besides rosette knots the following knots with low
crossing numbers are cylinder knots:
$9_{46}$, $9_{47}$, $10_{155}$, $3_1 \sharp 3_1^*$, $5_1 \sharp 5_1^*$,  
$7_1 \sharp 7_1^*$.   

\bigskip
{\bf Acknowledgements:} We thank D.~Zagier for his help proving Lemma 
\ref{independent_tan} and M.~Veve for the graphics in Figures
\ref{exeins} and \ref{zwanz}.

\end{document}